\documentclass[12pt]{article}
\usepackage{amsfonts}
\usepackage{amsthm}
\usepackage{url}
\usepackage{mathtools}
\usepackage[all]{xy}
\usepackage{amssymb}
\usepackage{color}
\usepackage{mathrsfs}
\usepackage{tikz}
\usepackage{hyperref}
\usetikzlibrary{calc,intersections,through,backgrounds}

\usepackage{pdfsync}
\newtheorem{thm}{Theorem}[section]
\newtheorem{prop}[thm]{Proposition}
\newtheorem{lemma}[thm]{Lemma}

\newtheorem{definition}[thm]{Definition}

\newcommand{\Z}{{\mathbb Z}} 
\newcommand{\Q}{{\mathbb Q}}

\newcommand{\Qp}{{\mathbb {Q}_{p}}}

\newcommand{\Sp}{{\mathrm{Spec} \:}}

\newcommand{\cI}{{\mathrm{c\text{--} Ind}}}
\DeclareMathOperator{\Spm}{\mathrm{m-Spec}}

\makeatletter
\let\c@equation\c@thm
\makeatother
\numberwithin{equation}{section}

\title{On the Breuil-Schneider conjecture II: Potentially crystalline non-generic case}

\author{Alexandre Pyvovarov}

\date{\today}

\begin{document}

\maketitle

\begin{abstract}

This paper improves some results of the author's previous work. We investigate the case of non-smooth points on automorphic components and prove the Breuil--Schneider conjecture under a generic Fontaine--Laffaille hypothesis on the residual Galois representation. The residual hypothesis is formulated in terms of an extremal modular Serre weight and the maximal shape of the associated Fontaine--Laffaille module; it is the hypothesis under which weight elimination makes the relevant patched intertwining operator invertible.

\end{abstract}

\tableofcontents

\section{Introduction}\label{Intro}
The aim of this work is to improve the results of \cite{Pyv3}, where the results depend on the assumption that a Galois representation comes from a point lying on an automorphic component. We will also use the patching construction of \cite{MR3529394}. The conjecture in question predicts the existence of $GL_n(F)$-invariant norms on locally algebraic $GL_n(F)$-representations, where $F$ is a $p$-adic field. This conjecture was first proposed by Breuil and Schneider in \cite{MR2359853}. For a brief survey of this conjecture the reader is referred to \cite{MR3409331}.

\subsection {Notation} \label{sec: Not}

Our setup is similar to \cite{Pyv3} and \cite{MR3529394}. Let $p$ a prime number such that $p\nmid 2n$. Let $F$ be a finite extension of $\Qp$ with a finite residue field $k_{F}$. Let $\mathcal{O}_F$ be its complete discrete valuation ring, let $\mathfrak{p}$ be the maximal ideal of $\mathcal{O}_F$ with uniformizer $\varpi$, and let $q=|\mathcal{O}_F/\varpi \mathcal{O}_F|$. Let $G=GL_n(F)$.

Let $E$ be a finite extension of $\Qp$ (the field of coefficients), $\mathcal{O}$ the ring of integers of $E$ and $\mathbb{F}$ the residue field. Fix a residual Galois representation  $\overline{r} : G_F \longrightarrow GL_n(\mathbb{F})$ of the local Galois group $G_F:= \mathrm{Gal}(\bar{F}/F)$. We assume that $E$ is large enough to contain all the embeddings $F \hookrightarrow \overline{\Q}_{p}$.

\subsection{Main result}\label{I.4}

Our setup is the following, we have fixed a residual Galois representation $\overline{r} : G_F \longrightarrow GL_n(\mathbb{F})$. The patching construction carried out in \cite{MR3529394} associates to $\overline{r}$ a modules $M_{\infty}(\overline{r}):=M_{\infty}$. It is enough to view $M_{\infty}$ as an $R_{\infty}[G]$-module, such that it is finitely generated over $R_{\infty}[[GL_n(\Z_p)]]$, where $R_\infty$(cf. \cite[section 2.8, p.27]{MR3529394}) is a complete Noetherian local $R_{\tilde{\mathfrak{p}}}^{\square}$-algebra with residue field $\mathbb{F}$. In \cite[section 5]{Pyv3} we define the rings $R_{\infty}(\sigma_{min})(\overline{r}):=R_{\infty}(\sigma_{min})$ and $R_{\infty}(\sigma_{max})(\overline{r}):=R_{\infty}(\sigma_{max})$.

Let $\rho : G_F \longrightarrow GL_n(E)$ be a potentially semi-stable Galois representation of weight $\sigma_{alg}$ and inertial type $\tau$. By the theory of Fontaine, the Galois representation $\rho$ corresponds to a filtered admissible $(\varphi,N, \mathrm{Gal}(L/F))$-module $\tilde{D}$. Proposition \ref{2.5} says that there is an admissible filtered $\varphi$-module $D$ such that the underlying Weil representation and the Hodge-Tate weights of $D$ are the same as that of $\tilde{D}$, but the monodromy operator is equal to zero. Let $r$ be the Galois representation corresponding to $D$. This Galois representation is potentially crystalline. Observe that the Galois representations $r$ and $\rho$ may have different reduction mod $p$ since the filtrations  on the corresponding $(\varphi,N)$-modules may differ. 

Let now $\overline{\rho} : G_F \longrightarrow GL_n(\mathbb{F})$ be the reduction mod $p$ of $\rho$. We can consider a new deformation problem associated to $\overline{\rho}$. Let $\tilde{M}_{\infty}:=M_{\infty}(\overline{\rho})$,   $\tilde{R}_{\infty}(\sigma_{min}):=R_{\infty}(\sigma_{min})(\overline{\rho})$ and let $\tilde{R}_{\infty}(\sigma_{min})'=R_{\infty}\otimes_{R_{\tilde{\mathfrak{p}}}^{\square}}R^{\square}_{\overline{\rho}}(\tau, \mathbf{v})$, where the ring $R^{\square}_{\overline{\rho}}(\tau, \mathbf{v})$ parametrizes all the potentially semi-stable lifts of $\overline{\rho}$ of weight $\sigma_{alg}$ and inertial type $\tau$.

The main result is the following theorem:

\begin{thm}\label{t1}
	Let $r$ and $\rho$ be two Galois representations, as above. Assume that $F/\Qp$ is unramified, that the inertial type $\tau$ is trivial, and that the Hodge type is associated with a $5n$-generic Fontaine--Laffaille weight. Assume moreover that the residual representation $\overline r$ satisfies the extremal Fontaine--Laffaille hypothesis of Definition \ref{def:GFL}. Assume that $\rho$ is generic (we do not assume that $\rho$ corresponds to a point lying on an automorphic component) and that $r_x:=r$ corresponds to a closed point $x \in \Sp R_{\infty}(\sigma_{max})[1/p]$ (deformation problem for $\overline{r}$). Then $BS(r_x)$ and $BS(\rho)$ both admit a $G$-invariant norm, where $BS(r_x)$ is a locally algebraic representation of $G$ encoding the Hodge--Tate weights of $r_x$ and the Weil--Deligne representation associated to $r_x$ (cf. \cite[section 2.3]{Pyv3}). The completions of $BS(r)$ and $BS(\rho)$ with respect to these norms are admissible.
\end{thm}

It is conjectured in \cite{MR3529394} that $V(r_x)$ depends only on the Galois representation $r_x$ and that $r_x \mapsto V(r_x)$ realizes the hypothetical $p$-adic local Langlands correspondence. Our Theorem \ref{t1} provides further evidence of this conjecture and generalizes the results of \cite{Pyv3}.

\subsection{Outline of the paper}

This article is organised as follows: In the section \ref{K.1}, we will recall definition of weakly admissible $(\varphi,N)$-modules with Galois action. Then in section \ref{K.2} we will prove that if we set $N=0$ on the weakly admissible $(\varphi, N)$ module $D$, then there is a filtration on the underlying $\varphi$-module of $D$, such that the $\varphi$-module then $D$ is again weakly admissible. The section \ref{A.2} is the heart of this paper. In section \ref{A.2.2} we will prove our main theorem and in section \ref{A.2.3} we will give an example when the Breuil-Schneider conjecture holds without assuming that a point lies on an automorphic component.

\section{Potentially crystalline representations}\label{K}

Let $D$ be  a weakly admissible $(\varphi, N)$-module. Here we will prove that if we set $N=0$ then there is a filtration on  $D$, the underlying $\varphi$-module of $D$, such that the $\varphi$-module then $D$ is again weakly admissible. First we recall few definitions and then we will study in detail the two dimensional case. Then we will see how some elementary inequalities allow us to deduce this result.

\subsection{Notation}\label{K.1}

Recall that $p$ is a prime number. In this section fix two finite extensions $F$ (the base field) and $E$ (the
coefficient field) of $\Qp$ such that $[F : \Qp] = |\mathrm{Hom}_{\Qp} (F, E)|$ where $\mathrm{Hom}_{\Qp}(F,E)$ denotes the set of all $\Qp$-linear embeddings of the field $F$ into the field $E$. We assume $F$ is contained in an algebraic closure $\overline{\Q}_{p}$ of $\Qp$. We denote by $q=p^{f_0}$ the cardinality of the residue field of $F$ and by $F_0 = \mathrm{Frac}(W(\mathbb{F}_q))$ its maximal unramified subfield. If $e := [L : \Qp]/f_0$, we set $val_F(x) := e.val_{\Qp}(x)$ (where $val_{\Qp}(p) := 1$) and $|x|_F := q^{-val_F(x)}$ for any $x$ in a finite extension of $\Qp$. We denote by $W_F=W(\overline{\Q}_{p}/F)$ (resp. $G_F:=\mathrm{Gal}(\overline{\Q}_{p}/F)$) the Weil (resp. Galois) group of $F$ and by
$\mathrm{rec}_p : W(\overline{\Q}_{p}/F)^{ab} \to F^{\times}$ the reciprocity map sending the geometric Frobenius to the uniformizer.

Let $L$ be a finite Galois extension of $L$ and $L_0$ its maximal unramified subfield. We assume $[L_0: \Qp] = |\mathrm{Hom}_{\Qp}(L_0,E)|$ and we let $p^{f}$ be the cardinality of the residue field of $L_0$ and $\varphi_0$ be the Frobenius on $F$ (raising to the $p$ each component of the Witt vectors). Consider the following two categories:

\begin{enumerate}
\item  the category $\mathrm{WD}_{L/F}$ of representations $(r, N, V)$ of the Weil-Deligne group of $F$ on a $E$-vector space $V$ of finite dimension such that $r$ is unramified when restricted to $W(\overline{\Q}_{p}/L)$.

\item the category $\mathrm{MOD}_{L/F}$ of quadruples $(\varphi, N, \mathrm{Gal}(L/F), D)$ where $D$ is a free $L_0 \otimes_{\Qp}E$-module of finite rank endowed with a Frobenius $\varphi : D \to D$, which is $\phi_0$-semi-linear bijective map, an $L_0 \otimes_{\Qp}E$-linear endomorphism $N : D \to D$ such that $N\varphi = p\varphi N$ and an action of $\mathrm{Gal}(L/F)$ commuting with $\varphi$ and $N$.
\end{enumerate}

There is a functor (due to Fontaine):
\[\mathrm{WD} : \mathrm{MOD}_{L/F} \longrightarrow \mathrm{WD}_{L/F}\]

The following proposition was proven in \cite{MR2359853}(Proposition 4.1):

\begin{prop}\label{2.1}
The functor $\mathrm{WD} : \mathrm{MOD}_{L/F} \longrightarrow \mathrm{WD}_{L/F}$ is an equivalence
of categories.
\end{prop} 

\noindent Denote $\mathrm{MOD}$ a quasi inverse of the functor $\mathrm{WD}$.

\medskip

\noindent If $D$ is an object of $\mathrm{MOD}_{L/F}$, we define:
\[ t_N(D)= \frac{1}{[F:L_0]f} val_F(\mathrm{det}_{L_0}(\varphi^{f}|D))\]

\noindent For $\sigma : L \hookrightarrow K$, let $D_{L}=D \otimes_{L_0}L$ and :
\[D_{L,\sigma}= D_{L} \otimes_{L\otimes_{\Qp}E} (L\otimes_{F,\sigma}E)\]

Then one has $D_{L} \simeq \prod_{\sigma: F \to E}D_{L,\sigma}$. To give an $L\otimes_{\Qp}E$-submodule $\mathrm{Fil}^{i}D_{L}$ of $D_{L}$ preserved by $\mathrm{Gal}(L/F)$ is the same thing as to give a collection $(\mathrm{Fil}^{i}D_{L,\sigma})_{\sigma}$ where $\mathrm{Fil}^{i}D_{L,\sigma}$ is a free $L\otimes_{F,\sigma}E$-submodule of $D_{L,\sigma}$ (hence a direct factor as $L\otimes_{F,\sigma}E$-module) preserved by the action of $\mathrm{Gal}(L/F)$. If $(\mathrm{Fil}^{i}D_{L,\sigma})_{\sigma,i}$ is a decreasing exhaustive separated filtration on $D_{L}$ by $L\otimes_{\Qp}E$-submodules indexed by $i \in\Z$ and preserved by $\mathrm{Gal}(L/F)$, we define:
\[t_H(D_{L}) = \sum_{\sigma} \sum_{i \in \Z} i \dim_{L}(\mathrm{Fil}^{i}D_{L,\sigma}/\mathrm{Fil}^{i+1}D_{L,\sigma})\]

Recall that such a filtration is called admissible if $t_H(D_{L})=t_N(D)$ and if, for any $L_0$-vector subspace $D'\subseteq D$ preserved by $\varphi$ and $N$ with the induced filtration on $D'_{L}$, one has $t_H(D'_{L})\leq t_N(D')$.

\subsection{Weakly admissible modules}\label{K.2}

\subsubsection{Examples}\label{K.2.1}

Let's consider the case of semi-stable representations of $\mathrm{Gal}(\overline{\Q}_{p}/\Qp)$. We know that the category of semi-stable representations is equivalent to the category of filtered weakly admissible ($\varphi$, $N$)-modules. In what follows we will deal exclusively with that type of ($\varphi$, $N$)-modules. 

Let's examine what happens in two dimensional case. In order to avoid specifying the extension $E$ over which the module $D$ is semi stable we will work with coefficients in $\overline{\Q}_{p}$. Assume that Hodge-Tate weights are $(0, k-1)$, with $k \geq 2$. The weakly admissible modules are of the form $D = \overline{\Q}_{p}e_{1}\oplus \overline{\Q}_{p}e_{2}$ with $\mathrm{Fil}^{i}D = D$ for $i \leq 0$, $\mathrm{Fil}^{1}D =\ldots=\mathrm{Fil}^{k-1} D \neq 0$ and $\mathrm{Fil}^{i}D = 0$ for $i \geq k$. Let $v_{p}$ denote the $p$-adic valuation, such that $v_{p}(p)=1$.

Consider the module given by:
$$
\left\{
    \begin{array}{lll}
        \varphi(e_{1}) & = & p \lambda e_{1}\\
        \varphi(e_{1}) & = & \lambda e_{2}\\
        \mathrm{Fil}^{k-1} D  & = & \overline{\Q}_{p}(e_{1}+\mathcal{L}e_{2})\\
        N(e_{1}) &=& e_{2}\\
        N(e_{2}) & =& 0 \\
        \mathcal{L} & \in & \overline{\Q}_{p}\\
    \end{array}
\right.
$$

The non-trivial only $(\varphi, N)$-stable sub vector spaces is $\overline{\Q}_{p}e_{2}$ so the admissibility condition is equivalent to $k-1=2v_p(\lambda)+1$ and $0=t_H(\overline{\Q}_{p}e_{2})\leq t_N(\overline{\Q}_{p}e_{2})=v_p(\lambda)$.

If we set $N=0$ and we want to keep the same $\varphi$ and the same filtration, then $D$ the form:
$$
\left\{
    \begin{array}{lll}
        \varphi(e_{1}) & = & p \lambda e_{1}\\
        \varphi(e_{1}) & = & \lambda e_{2}\\
        \mathrm{Fil}^{k-1} D  & = & \overline{\Q}_{p}(e_{1}+e_{2})\\

    \end{array}
\right.
$$

Notice this is the same filtration as above but we made a base change $e_{2}\mapsto \mathcal{L} e_{2}$. From the admissibility of $(\varphi,N)$-module we get the admissibility of $\varphi$-module. Indeed  the non-trivial $\varphi$-stable sub vector spaces are $\overline{\Q}_{p}e_{1}$, $\overline{\Q}_{p}e_{2}$ and we have $0=t_H(\overline{\Q}_{p}e_{1}) < t_N(\overline{\Q}_{p}e_{1})=v_p(\lambda)+1$. This was the only inequality that was left to check.

For 3-dimensional examples the reader may look at \cite{MR3459578} and check in the same way that for every isomorphism class with non zero monodromy operator, we get a crystalline representation when we kill the monodromy operator.

\subsubsection{Inequalities of integers}\label{K.2.2}

\begin{lemma}\label{2.2}
Let $i_1 \leq i_2 \leq \ldots \leq i_d$, a sequence of integers in $\Z$ and let $c\in \Q$. Assume that we have
$$ \sum\limits_{j=1}^{d} i_j \leq dc$$

\noindent then for all $1\leq n \leq d$, we have
$$ \sum\limits_{j=1}^{n} i_j \leq nc$$
\end{lemma}

\begin{proof} We have that $i_{n+1} \geq i_j$ for all $1\leq j \leq n$. Then it follows that $ni_{n+1} \geq \sum_{j=1}^{n} i_j$. By a simple calculation this inequality is equivalent to 
$$\frac{1}{n+1}\sum_{j=1}^{n+1} i_j - \frac{1}{n}\sum_{j=1}^{n} i_j =  \frac{i_{n+1}}{n+1}+(\frac{1}{n+1}-\frac{1}{n}) \sum_{j=1}^{n} i_j$$
$$=\frac{1}{n(n+1)}(ni_{n+1} - \sum_{j=1}^{n} i_j)\geq 0$$

\noindent Then it follows by induction that 
$$\frac{1}{n}\sum_{j=1}^{n} i_j \leq \frac{1}{d}\sum_{j=1}^{d} i_j \leq c$$

\noindent The result follows.
\end{proof}

\begin{lemma}\label{2.3}
Let $i_1 \leq i_2 \leq \ldots \leq i_{d+k}$, a sequence of integers in $\Z$ and let $c_{1}$, $c_{2}\in \Q$ . Assume that we have
$$ \sum\limits_{j=1}^{d} i_j \leq dc_{1}$$

\noindent and
$$ \sum\limits_{j=1}^{d} i_j + \sum\limits_{j=d+1}^{d+k} i_j \leq dc_{1} +k(c_{2}+1)$$

\noindent then for all $1\leq n \leq k$, we have
$$ \sum\limits_{j=1}^{d} i_j + \sum\limits_{j=d+1}^{d+n} i_j\leq  dc_{1} +n(c_{2}+1)$$
\end{lemma}

\begin{proof} We prove first the desired inequality for $n=1$. If $i_{d+1} \leq c_{2}+1$ then 
$$ \sum\limits_{j=1}^{d} i_j + i_{d+1}\leq dc_{1} + i_{d+1} \leq  dc_{1} +(c_{2}+1)$$

\noindent If $i_{d+1} \geq c_{2}+1$, then 
$$ \sum\limits_{j=1}^{d} i_j + k i_{d+1} \leq \sum\limits_{j=1}^{d} i_j + \sum\limits_{j=d+1}^{d+k} i_j \leq dc_{1} +k(c_{2}+1)$$

Subtracting  $(k-1)i_{d+1} + c_{2}+1 + dc_{1}$ on both sides of the previous inequality, we get 
$$ \sum\limits_{j=1}^{d} i_j +  i_{d+1} -(dc_{1}+c_{2}+1) \leq (k-1)(-i_{d+1}+c_{2}+1) \leq 0$$

\noindent Hence in any case we have that 
$$ \sum\limits_{j=1}^{d} i_j + i_{d+1} \leq  dc_{1} +(c_{2}+1)=d(c_{1} + (c_{2}+1)/d)$$

This proves the lemma for $n=1$. Now replacing $d$ by $d+1$, $c_{1}$ by $c_{1} + (c_{2}+1)/d$ and $k$ by $k-1$, so that we can repeat the procedure above for $n=1$, i.e. we start now with inequalities:
$$ \sum\limits_{j=1}^{d+1} i_j \leq d(c_{1} + (c_{2}+1)/d)$$

\noindent and
$$ \sum\limits_{j=1}^{d+1} i_j + \sum\limits_{j=d+2}^{d+k} i_j \leq d(c_{1} + (c_{2}+1)/d) +(k-1)(c_{2}+1)$$

\noindent then proceed as before to get:
$$ (\sum\limits_{j=1}^{d} i_j + i_{d+1}) + i_{d+2} \leq d(c_{1} + (c_{2}+1)/d) +(c_{2}+1)= dc_{1} +2(c_{2}+1),$$

\noindent i.e. an inequality for $n=2$. We proceed by induction in a similar fashion to prove this lemma.
\end{proof}

By induction from the two previous lemmas we get the following result:

\begin{lemma}\label{2.4}
Let $i_1 \leq i_2 \leq \ldots \leq i_{d_{1}+\ldots+d_{s}}$, a sequence of integers in $\Z$ and let $c_{i}\in \Q$ for all $1\leq i \leq s$. Assume that we have
$$ \sum\limits_{j=1}^{d_{1}} i_j \leq d_{1}c_{1}$$
$$ \sum\limits_{j=1}^{d_{1}} i_j + \sum\limits_{j=d_{1}+1}^{d_{1}+d_{2}} i_j \leq d_{1}c_{1} +d_{2}c_{2}$$
$$\vdots$$
$$ \sum\limits_{j=1}^{d_{1}} i_j +\ldots+ \sum\limits_{j=d_{1}+\ldots d_{s-1}+1}^{d_{1}+\ldots+d_{s}} i_j \leq d_{1}c_{1}+\ldots +d_{s}c_{s}$$

\noindent then for all $1\leq k \leq s$ and for all $1\leq n \leq d_{k}$, we have
$$ \sum\limits_{j=1}^{d_{1}} i_j +\ldots +\sum\limits_{j=d_{1}+\ldots d_{k-1}+1}^{d_{1}+\ldots d_{k-1}+n} i_j\leq  d_{1}c_{1}+\ldots+d_{k-1}c_{k-1} +nc_{k}$$
\end{lemma}

\subsubsection{General case}\label{K.2.3}

\begin{prop}\label{2.5}
Let $(\varphi, N, \mathrm{Gal}(L/F), D)$ be an object in the category $\mathrm{MOD}_{L/F}$ which has an admissible filtration. Then the object
$(\varphi, 0, \mathrm{Gal}(L/F), D)$ has also an admissible filtration. The Hodge-Tate weights of $(\varphi, 0, \mathrm{Gal}(L/F), D)$ and  $(\varphi, N, \mathrm{Gal}(L/F), D)$ are the same and those objects have the same action of $\varphi$ on $D$.
\end{prop}

\begin{proof} Let $(r,N,V)= \mathrm{WD}(D)$. Assume that the extension $E$ of $\Qp$ is big enough so that 
\[(r,N,V)=\bigoplus_{i=1}^{s} (r_i,N_i,V_i)\]

\noindent where $(r_i, N_i, V_i)$ is absolutely indecomposable of dimension $d_i$. 

Let $(\varphi_i,N_i, D_i) =\mathrm{MOD}((r_i,N_i,V_i))$, then it is an absolutely indecomposable object in $\mathrm{MOD}_{L/F}$. Let $D_{i,0} = \mathrm{Ker}(N_i : D_i \rightarrow D_i)$ and $\varphi_{i,0}=\varphi_{i}|D_{i,0}$. An indecomposable Weil-Deligne representation can be always written in a specific form (cf. 3.1.3 (ii) \cite{MR0389771}). Then $D_i$ can be written as
\[D_i=D_{i,0}\oplus D_{i,0}(1) \oplus \ldots \oplus D_{i,0}(b_i-1)\]

\noindent where $D_{i,0}(n)\simeq D_{i,0}$ with $\varphi_{i,0}|D_{i,0}(n)=p^{n}\varphi_{i,0}$, $N_{i}|D_{i,0} = 0 $ and $N_i$ sends $D_{i,0}(n)$ to $D_{i,0}(n-1)$ via identity if $n>0$. Note that $D_{i,0}$ is absolutely irreducible. Since $\varphi_{i,0}^{f}$ is $L_{0}\otimes_{\Qp}E$-linear and commutes with $\mathrm{Gal}(L/F)$ $\varphi_{i,0}$, then $\varphi_{i,0}^{f}$ is a scalar matrix with values in $F^{\times}$. Let $n_i
= \dim D_{i,0}$ and $\varphi_{i,0}^{f}=\lambda_i.\mathrm{Id}$.

Choose an order on summands such that $\mathrm{val}_{F}(\lambda_1) \leq \mathrm{val}_{F}(\lambda_2) \leq \ldots \leq \mathrm{val}_{F}(\lambda_s)$.

For an embedding $\sigma$ the Hodge-Tate weights are $i_{\sigma,1}<\ldots<i_{\sigma,n}$. Write $c_i=\frac{1}{[E:L_0]f}\mathrm{val}_{F}(\lambda_i)$ and $i_j = \sum_{\sigma} i_{\sigma, j}$. 

Notice that the only sub-objects of $D_i$ are $D_{i,0} \oplus \ldots \oplus D_{i,0}(k)$ for $0\leq k \leq b_i-1$. Then the admissibility condition of $D$, for these subobjects, gives us the following inequalities:

\begin{enumerate}
\item[1.] admissibility for $D_{1,0} \oplus \ldots \oplus D_{1,0}(k)$
\[\sum\limits_{j=1}^{n_1} i_j + \ldots +\sum\limits_{j=kn_1+1}^{(k+1)n_1}i_j \leq n_1 c_1+ \ldots + n_1(c_1+k)\]
\noindent for $0\leq k \leq b_1-1$,

\item[2.] admissibility for $D_1 \oplus D_{2,0} \oplus \ldots \oplus D_{2,0}(k)$
\[\sum\limits_{j=1}^{d_1} i_j + \sum\limits_{j=d_1+1}^{d_1+n_2} i_j + \ldots +\sum\limits_{j=d_1+kn_2+1}^{d_1+(k+1)n_2}i_j \leq\] \[\leq n_1 c_1+ \ldots + n_1(c_1+b_1-1)+n_2 c_2+ \ldots + n_2(c_2+k)\]
\noindent for $0\leq k \leq b_2-1$,

$$\vdots$$

\item[s.] admissibility for $D_1 \oplus \ldots \oplus D_{s-1} \oplus D_{s,0} \oplus \ldots \oplus D_{s,0}(k)$
\[\sum\limits_{j=1}^{d_1+\ldots+d_{s-1}} i_j + \sum\limits_{j=d_1+\ldots+d_{s-1}+1}^{d_1+\ldots+d_{s-1}+n_s} i_j + \ldots +\sum\limits_{j=d_1+\ldots+d_{s-1}+kn_s+1}^{d_1+\ldots+d_{s-1}+(k+1)n_s}i_j \leq\] \[\leq \sum\limits_{i=1}^{s-1} \sum\limits_{l=0}^{b_i-1} n_i (c_i+l)+ n_s c_s+ \ldots + n_s(c_s+k)\]
\noindent for $0\leq k \leq b_s-1$, with an equality for $k=b_s$. 
\end{enumerate}

\noindent Then applying Lemma \ref{2.4} to each set of inequalities above, from 1 to s, we get the following intermediate inequalities:

\begin{enumerate}
\item[1.] For $0\leq k \leq b_1-1$,
\[\sum\limits_{j=1}^{a} i_j  \leq a c_1\]
\noindent for $1\leq a \leq n_1$,
\[\vdots\]
\[\sum\limits_{j=1}^{n_1} i_j + \ldots +\sum\limits_{j=kn_1+1}^{a}i_j \leq n_1 c_1+ \ldots + a(c_1+k)\]
\noindent for $kn_1+1\leq a \leq (k+1)n_1$,

\item[2.] For $0\leq k \leq b_2-1$,
\[\sum\limits_{j=1}^{d_1} i_j + \sum\limits_{j=d_1+1}^{d_1+a} i_j  \leq n_1 c_1+ \ldots + n_1(c_1+b_1-1)+a c_2\]
\noindent for $d_1+1\leq a \leq n_2$,
\[\vdots\]
\[\sum\limits_{j=1}^{d_1} i_j + \sum\limits_{j=d_1+1}^{d_1+n_2} i_j + \ldots +\sum\limits_{j=d_1+kn_2+1}^{a}i_j \leq\] \[ \leq n_1 c_1+ \ldots + n_1(c_1+b_1-1)+n_2 c_2+ \ldots + a(c_2+k)\]
\noindent for $d_1+kn_2+1\leq a \leq d_1+(k+1)n_2$,

\[\vdots\]

\item[s.] For $0\leq k \leq b_s-1$,
\[\sum\limits_{j=1}^{d_1+\ldots+d_{s-1}} i_j + \sum\limits_{j=d_1+\ldots+d_{s-1}+1}^{a} i_j  \leq \sum\limits_{i=1}^{s-1} \sum\limits_{l=0}^{b_i-1} n_i (c_i+l)+ a c_s\]
\noindent for $d_1+\ldots+d_{s-1}+1\leq a \leq d_1+\ldots+d_{s-1}+n_s$,
\[\vdots\]
\[\sum\limits_{j=1}^{d_1+\ldots+d_{s-1}} i_j + \sum\limits_{j=d_1+\ldots+d_{s-1}+1}^{d_1+\ldots+d_{s-1}+n_s} i_j + \ldots +\sum\limits_{j=d_1+\ldots+d_{s-1}+kn_s+1}^{d_1+\ldots+d_{s-1}+a}i_j \leq\] \[\leq \sum\limits_{i=1}^{s-1} \sum\limits_{l=0}^{b_i-1} n_i (c_i+l)+ n_s c_s+ \ldots + a(c_s+k)\]
\noindent for $d_1+\ldots+d_{s-1}+kn_s+1\leq a \leq d_1+\ldots+d_{s-1}+(k+1)n_s$.
\end{enumerate}

By equivalence $(i)\Leftrightarrow (ii)$ of Proposition 3.2 \cite{MR2359853}, all these inequalities tell that there is an admissible filtration on the $\varphi$-modules $(\varphi,D^{N=0})$. By construction the $\varphi$-modules $(\varphi,D^{N=0})$ has the same Hodge-Tate weights as $(\varphi,N)$-module $D$ and both modules inherit the same action of $\varphi$.
\end{proof}

\section{Potentially crystalline non-generic points}\label{A.2}

In this section we will prove the existence of a $G$-invariant norm on $BS(r)$ in some cases when $r$ is a potentially crystalline Galois representation which is not necessarily generic. A more precise statement will be given in Theorem \ref{3.4}. Similarly to \cite{Pyv3} we will embed $BS(r)$ into a unitary $E$-Banach space representation of $G$. In section \ref{A.2.1} we build a framework for the proof by examining the support of the patched modules $M_{\infty}(\sigma_{min}^{\circ})$ and $M_{\infty}(\sigma_{max}^{\circ})$. The additional input is a generic Fontaine--Laffaille weight-elimination theorem for the residual representation. Then in section \ref{A.2.2} we prove the main result. In the last section we give an example illustrating Theorem \ref{3.4}. We refer the reader to \cite[section 5]{Pyv3} for the definition of $M_{\infty}(\sigma_{min}^{\circ})$ and $M_{\infty}(\sigma_{max}^{\circ})$.

\subsection{More on support of patched modules}\label{A.2.1}

Let $x \in \Spm R_{\infty}(\sigma_{min})[1/p]$ and $y$ be the image of $x$ in $\Sp \mathfrak{Z}_{\Omega}$ by the map $\alpha^{\sharp}$ from \cite[Theorem 5.6]{Pyv3}. Define 
$$\gamma_x:=\cI_K^G \sigma_{max}(\lambda) \otimes_{\mathfrak{Z}_{\Omega}} \kappa(y)$$

\noindent and 
$$\delta_x:= \cI_K^G \sigma_{min}(\lambda) \otimes_{\mathfrak{Z}_{\Omega}} \kappa(y).$$

Let now $x \in \Spm R_{\infty}(\sigma_{max})[1/p]$. To the point $x$ corresponds a Galois representation $r_x$. Let $\Sigma$ be a subset of $\Spm R_{\infty}(\sigma_{max})[1/p]$, consisting of those $x$ such that the representation $\pi_{sm}(r_x)$ is generic. We will prove that the set $\Sigma$ is a dense subset of $\Spm R_{\infty}(\sigma_{max})[1/p]$.

\noindent Let $H:=\mathrm{Hom}_{G}(\mathrm{c\text{--} Ind}_{K}^{G} \sigma_{min}(\lambda), \mathrm{c\text{--} Ind}_{K}^{G} \sigma_{max}(\lambda))$. Consider the natural evaluation map:
$$ev  :  H \otimes_{\mathfrak{Z}_{\Omega}} \mathrm{c\text{--} Ind}_{K}^{G} \sigma_{min}(\lambda)\rightarrow \mathrm{c\text{--} Ind}_{K}^{G} \sigma_{max}(\lambda),$$

\noindent given by $f \otimes v  \mapsto  f(v)$. It follows from \cite[Theorem 7.1]{Pyv1}, that $H$ is locally free $\mathfrak{Z}_{\Omega}$-module of rank one.

Recall that $\sigma_{max}:= \sigma_{max}(\lambda) \otimes \sigma_{alg}$ and $\sigma_{min}:= \sigma_{min}(\lambda) \otimes \sigma_{alg}$. We have an isomorphism:
$$H \simeq \mathrm{Hom}_G(\cI_K^{G} \sigma_{min}, \cI_K^{G} \sigma_{max}),$$

\noindent thus $\mathrm{Hom}_G(\cI_K^{G} \sigma_{min}, \cI_K^{G} \sigma_{max})$ is a locally free $\mathfrak{Z}_{\Omega}$-module of rank one. Let $\phi$ be the image of $ev$ by the functor $\mathrm{Hom}_{E}^{cont}\left(\mathrm{Hom}_{G}(., (M_{\infty})^{d}[1/p]),E\right)$. Then:
$$\phi : \mathrm{Hom}_G(\cI_K^{G} \sigma_{min}, \cI_K^{G} \sigma_{max}) \otimes_{\mathfrak{Z}_{\Omega}} M_{\infty}(\sigma_{min}^{\circ})[1/p] \longrightarrow M_{\infty}(\sigma_{max}^{\circ})[1/p]$$

\noindent is a homomorphism of $R_{\infty}(\sigma_{min})$-modules. Let $r_x$ be the Galois representations corresponding to the point $x$. Then by Proposition 4.33 \cite{MR3529394}, we have $$V(r_x)^{l.alg}=\pi_x\otimes \pi_{alg}(r_x),$$ \noindent where $\pi_x$ is some smooth admissible representation in the Bernstein component $\Omega$.

Let $X$ be the set of points $x$ such that $\phi \otimes \kappa(x) \neq 0$.

\begin{lemma}\label{3.2}
By assumption we have that $x\in \mathrm{Supp}(M_{\infty}(\sigma_{max}^{\circ}))$. It follows $\mathrm{Hom}_K(\sigma_{max},V(r_x)^{l.alg})\neq 0$ so that we have a non-zero map $\gamma_x \rightarrow \pi_x$. Then $x\in X$ if and only the composition $\delta_x \xrightarrow{\Delta} \gamma_x \rightarrow \pi_x$ is non-zero for some (equivalently any) non-zero $\Delta \in \mathrm{Hom}_G(\delta_x, \gamma_x)$.
\end{lemma}

\begin{proof} If $x\in X$, then the specialization 
$$\phi \otimes \kappa(x) : M_{\infty}(\sigma_{min}^{\circ}) \otimes_{R_{\infty}} \kappa(x) \rightarrow M_{\infty}(\sigma_{max}^{\circ}) \otimes_{R_{\infty}} \kappa(x)$$ 

\noindent is non zero, where $\kappa(x)$ is the residue field at $x$. However by Proposition 2.22 \cite{MR3306557} and Frobenius reciprocity we have:
$$M_{\infty}(\sigma_{min}^{\circ}) \otimes_{R_{\infty}} \kappa(x) \simeq \mathrm{Hom}_{E}^{cont}\left(\mathrm{Hom}_{G} (\mathrm{c\text{--} Ind}_{K}^{G} \sigma_{min}\otimes_{\mathfrak{Z}_{\Omega}}\kappa(y), V(r_{x})^{l.alg}),E\right)\simeq$$
$$\mathrm{Hom}_{E}^{cont}\left(\mathrm{Hom}_{G} (\delta_x \otimes \pi_{alg}(r_x), \pi_x\otimes \pi_{alg}(r_x)),E\right)\simeq\left(\mathrm{Hom}_{G} (\delta_x, \pi_x) \right)^{\ast}$$
\noindent where $(.)^{\ast} = \mathrm{Hom}_{E}(.,E)$ is the dual of finite dimensional vector spaces and similarly, $$M_{\infty}(\sigma_{max}^{\circ}) \otimes_{R_{\infty}} \kappa(x) \simeq \left(\mathrm{Hom}_{G} (\gamma_x, \pi_x) \right)^{\ast}.$$ 

\noindent It follows that the map $\phi \otimes \kappa(x)$ is induced by the following map:
$$\begin{array}{ccccc}
 &  & \mathrm{Hom}_G(\delta_x, \gamma_x) \otimes \mathrm{Hom}_{G}(\gamma_x, \pi_x) & \to & \mathrm{Hom}_{G}(\delta_x, \pi_x) \\
 & & \Delta \otimes f &  \mapsto & f \circ \Delta \\
\end{array}$$ 

\noindent The assertion of this lemma follows. Since $\mathrm{Hom}_G(\delta_x, \gamma_x)$ is one dimensional, any non-zero $\Delta \in \mathrm{Hom}_G(\delta_x, \gamma_x)$ will do.
\end{proof}

\begin{lemma} \label{3.3}

Let $x\in \mathrm{Supp}(M_{\infty}(\sigma_{max}^{\circ}))$ be a closed point. The following assertions are equivalent:
\begin{enumerate}
\item $x\in X$
\item The $G$-equivariant map $\gamma_x \rightarrow \pi_x$ is injective.
\end{enumerate}
\end{lemma}

\begin{proof} \underline{\textit{1} implies \textit{2}}. First notice that $\phi \otimes \kappa(x) \neq 0 \Longrightarrow ev \otimes \kappa(y) \neq 0$. Let $\mathrm{soc}_{G}(\gamma_x)$ be the  $G$-socle of $\gamma_x$ and let $\iota$ the image of the map $ev \otimes \kappa(y)$. The image $\iota$ has a finite length because the representation $\gamma_x$ is of finite length. Let $\nu$ be an irreducible quotient of $\iota$. Then  $\mathrm{Hom}_{G}(\delta_x, \nu) = \mathrm{Hom}_{G}(\mathrm{c\text{--} Ind}_{K}^{G} \sigma_{min}(\lambda), \nu) =\mathrm{Hom}_{K}(\sigma_{min}(\lambda), \nu) \neq 0$ because $ev \otimes \kappa(y) \neq 0$. By  \cite[Theorem 2.1]{Pyv2} the representation $\nu$ is generic, but Corollary 3.11 in \cite{MR3529394} says that the only irreducible generic subquotient of $\gamma_x$ is $\mathrm{soc}_{G}(\gamma_x)$, so $\nu = \mathrm{soc}_{G}(\gamma_x)$.

Then the map $ev \otimes \kappa(y)$ factors through $\mathrm{soc}_{G}(\gamma_x)$ so that the diagram below commutes:
$$\xymatrix{
\delta_x \ar[r]  \ar@{->>}[d] \ar@{->>}[rd]
&\gamma_x \\
\iota \ar@{->>}[r] &\mathrm{soc}_{G}(\gamma_x) \ar@{^{(}->}[u]
}$$
In particular the composition:
$$\delta_x \twoheadrightarrow \mathrm{soc}_{G}(\gamma_x)\hookrightarrow \gamma_x \longrightarrow \pi_x$$
\noindent is non zero, by Lemma \ref{3.2}.

If the map $\gamma_x \longrightarrow \pi_x$ is not injective, let $\kappa$ be it's kernel. Since $\kappa$ is non zero by assumption it is equal or contains an irreducible representation $\eta$. The representation $\eta$ is also a sub-representation of $\gamma_x$. Since $\mathrm{soc}_{G}(\gamma_x)$ is irreducible it is a unique irreducible sub-representation of $\gamma_x$, hence $\eta = \mathrm{soc}_{G}(\gamma_x)$. Therefore $\mathrm{soc}_{G}(\gamma_x) \subseteq \kappa$, so the image of $\mathrm{soc}_{G}(\gamma_x)$ by the map $\gamma_x \longrightarrow \pi_x$ is 0. Since $\mathrm{soc}_{G}(\gamma_x)$ is irreducible the composite map $\mathrm{soc}_{G}(\gamma_x)\hookrightarrow \gamma_x \longrightarrow \pi_x$ is injective and 0.

This would imply that the composition:
$$\delta_x \twoheadrightarrow \mathrm{soc}_{G}(\gamma_x)\hookrightarrow \gamma_x \longrightarrow \pi_x$$
is  $0$, which is a contradiction. Therefore assertion \textit{2} follows.

\underline{\textit{2} implies \textit{1}}. Since $\mathrm{soc}_{G}(\gamma_x)$ is generic, by \cite[Theorem 2.1]{Pyv2} we have that $\mathrm{Hom}_G
(\delta_x, \mathrm{soc}_{G}(\gamma_x))\neq 0$. Then the composition $\delta_x \twoheadrightarrow \mathrm{soc}_{G}(\gamma_x)\hookrightarrow \gamma_x \hookrightarrow \pi_x$ is non zero and by Lemma \ref{3.2}, $x\in X$.

\end{proof}

Recall that $\Sigma$ is a subset of $\Spm R_{\infty}(\sigma_{max})[1/p]$, consisting of those $x$ such that the representation $\pi_{sm}(r_x)$ is generic.

\begin{lemma}\label{3.5}
We have the following inclusion $\Sigma \subseteq X$.
\end{lemma}

\begin{proof} Let $x \in \Sigma$. Since $x \in \mathrm{Supp}(M_{\infty}(\sigma_{max}^{\circ}))$, we have $\mathrm{Hom}_G(\gamma_x,\pi_x) \neq 0$. However, by Corollary 3.12 \cite{MR3529394}, we have $\pi_{sm}(r_x) \simeq \gamma_x$. It follows that $\gamma_x$ is irreducible. Thus any non-zero $G$-equivariant map $\gamma_x \rightarrow \pi_x$ is injective. The assertion follows from the Lemma \ref{3.3}.
\end{proof}

\begin{lemma}\label{3.7}
Let $S$ be an equidimensional Noetherian ring, and $\mathfrak{a}$ an ideal in $S$. Assume that $S$ is Jacobson. Then we have the following assertions
\begin{enumerate}
\item $\dim (S/\mathfrak{a}) = \dim S$ if and only if $V(\mathfrak{a})\cap \Spm S$ contains an irreducible component of $\Spm S$.
\item $\dim S/\mathfrak{a} < \dim S$ if and only if $\Spm S\setminus V(\mathfrak{a})\cap \Spm S$ is Zariski dense in $\Spm S$.
\end{enumerate}
\end{lemma}
\begin{proof} The assertions 1 and 2 are trivially equivalent. Let's prove the first one.

Assume that $\dim (S/\mathfrak{a}) = \dim S$ and write $V(\mathfrak{a})\cap \Spm S=\bigcup_{i} V(\mathfrak{q}_i)$ as a union of irreducible components. Then there is an index $i$ such that $\dim V(\mathfrak{q}_i) = \dim S$, it follows that $\mathfrak{q}_i$ is actually a minimal prime in $S$. Then $V(\mathfrak{q}_i)$ is an irreducible component of $\Spm S$.

Assume now that $V(\mathfrak{a})\cap \Spm S$ contains an irreducible component $V(\mathfrak{p})$ of $\Spm S$. From inclusions  $ V(\mathfrak{p}) \subseteq V(\mathfrak{a})\cap \Spm S \subseteq \Spm S$, follows that $ \dim  V(\mathfrak{p}) \leq \dim V(\mathfrak{a})\cap \Spm S \leq \dim S$. Since $S$ is equidimensional $\dim  V(\mathfrak{p}) = \dim S$, it follows that $\dim (S/\mathfrak{a})$  $=\dim S$.
\end{proof}

\begin{lemma}\label{3.8}
The set $\{x \in \Spm R_{\tilde{\mathfrak{p}}}^{\square}(\sigma_{max})[1/p] \:|\: \pi_{sm}(r_x) \mbox{ is generic} \}$ is dense in $ \Spm R_{\tilde{\mathfrak{p}}}^{\square}(\sigma_{max})[1/p]$.
\end{lemma}

\begin{proof} By Theorem 1.2.7 \cite{MR3546966} $x \in \Spm  R_{\tilde{\mathfrak{p}}}^{\square}(\sigma_{min})[1/p]$ is smooth if and only if $WD(r_x)$ is generic and by Lemma 1.1.3 \cite{MR3546966}, $WD(r_x)$ is generic if and only if $\pi_{sm}(r_x)$ is generic. Let $\mathcal{S}$ the singular locus of $\Sp R_{\tilde{\mathfrak{p}}}^{\square}(\sigma_{min})[1/p]$. Then we have that:
\[\{x \in \Spm R_{\tilde{\mathfrak{p}}}^{\square}(\sigma_{max})[1/p] \:|\: \pi_{sm}(r_x) \mbox{ is generic} \}\]
\[= \Spm R_{\tilde{\mathfrak{p}}}^{\square}(\sigma_{max})[1/p] \setminus (\Spm R_{\tilde{\mathfrak{p}}}^{\square}(\sigma_{max})[1/p] \cap \mathcal{S})\]

We know that the ring $ R_{\tilde{\mathfrak{p}}}^{\square}(\sigma_{min})[1/p]$ is equidimensional and the complement of the singular locus is Zariski dense, by Theorem (3.3.4) \cite{MR2373358}. Then by previous lemma we have that $\dim \mathcal{S} < \dim R_{\tilde{\mathfrak{p}}}^{\square}(\sigma_{min})[1/p]$. Moreover we have 
$$\dim (\Spm R_{\tilde{\mathfrak{p}}}^{\square}(\sigma_{max})[1/p] \cap \mathcal{S}) \leq \dim \mathcal{S} < \dim R_{\tilde{\mathfrak{p}}}^{\square}(\sigma_{min})[1/p]$$

\noindent and $\dim R_{\tilde{\mathfrak{p}}}^{\square}(\sigma_{min})[1/p] = \dim R_{\tilde{\mathfrak{p}}}^{\square}(\sigma_{max})[1/p]$. We conclude by Lemma \ref{3.7}, because the ring $R_{\tilde{\mathfrak{p}}}^{\square}(\sigma_{max})[1/p]$ is also equidimensional by Theorem (3.3.8) \cite{MR2373358}.
\end{proof}

\begin{prop}\label{3.9}
$\Sigma$ is Zariski dense in $\Sp R_{\infty}(\sigma_{max})[1/p]$. 
\end{prop}

\begin{proof} To prove that $\Sigma$ is Zariski dense in $\mathrm{Spec}(R_{\infty}(\sigma_{\max})[1/p])$, it would suffice to prove that those points are dense in every irreducible component of $R_{\infty}(\sigma_{max})[1/p]$. Since the spectrum of this ring is a union of irreducible components of the spectrum of $R_{\infty}(\sigma_{\max})'[1/p]$ by \cite[Lemma 4.18]{MR3529394}, it is enough to prove it for this ring. The result follows from previous lemma and the fact that closed points are dense in a Jacobson ring.
\end{proof}

We will prove now that the set $X$ is Zariski closed. We start with the following lemmas:

\begin{lemma}\label{3.12}
	The set $X$ is closed if and only if $\phi$ is surjective. In this case $X = \Sp R_{\infty}(\sigma_{max})[1/p]$.
\end{lemma}

\begin{proof} If $\phi$ is surjective then $X = \Sp R_{\infty}(\sigma_{max})[1/p]$. 
	
	If $X$ is closed then by Lemma \ref{3.5}, it contains the closure of $\Sigma$. However by Proposition \ref{3.9}, $X = \Sp R_{\infty}(\sigma_{max})[1/p]$. By \cite[Lemma 4.18(1)]{MR3529394}, $M_\infty(\sigma_{max}^{\circ})[1/p]$ is locally free of rank one on its support. Thus $\phi\otimes\kappa(x)$ is surjective at every closed point $x$. The cokernel of $\phi$ is finite; if it were non-zero, its support would contain a closed point. Hence the cokernel is zero and $\phi$ is surjective.
\end{proof}

\paragraph{A Fontaine--Laffaille input.}
From now until the end of this subsection assume that $F/\Qp$ is unramified and that the inertial type is trivial. Put $K=GL_n(\mathcal O_F)$ and $\eta=(n-1,\ldots,1,0)$. Let $F(\lambda)$ be the Serre weight associated with the Hodge type. For every embedding $j:k_F\hookrightarrow\mathbb F$, write $\lambda_j=(\lambda_{j,1},\ldots,\lambda_{j,n})$.

We retain the Taylor--Wiles and adequacy hypotheses used in \cite{MR3529394} to construct $M_\infty$; in particular, its action is arithmetic and $M_\infty$ is projective over $\mathcal O[[K]]$. After base change from $\mathfrak Z_\Omega$ to the completed local patched deformation ring, choose saturated lattices in the two types and a primitive integral normalization of the evaluation map. Write $\mathfrak Z_\Omega^\circ$ for the resulting completed integral Bernstein centre and $H^\circ$ for the saturated rank-one lattice in $H$. The integral evaluation induces a homomorphism between finite modules
\[
 \phi^\circ:H^\circ\widehat\otimes_{\mathfrak Z_\Omega^\circ}
 M_\infty(\sigma_{min}^\circ)\longrightarrow
 M_\infty(\sigma_{max}^\circ)
\]
whose generic fibre is $\phi$. Via the equivalence defined by the type, a reduced expression for the affine Weyl group element carrying the minimal type to the maximal type gives the usual factorization of $\phi^\circ$ into normalized rank-one intertwining operators; compare \cite[Proposition 2.2]{MR782228} and \cite[Theorem 7.6.20]{MR1204652}. We use the integral rank-one diagrams and normalizations of \cite[Section 10.2]{LLHMPQ}. Denote the $F(\lambda)$-relevant tame types occurring in these diagrams by $\tau_1,\ldots,\tau_m$. Multiplying the primitive normalization by a unit does not affect surjectivity.

\begin{definition}\label{def:GFL}
We say that $\overline r$ satisfies the \emph{extremal generic Fontaine--Laffaille hypothesis} (for the above evaluation map) if the following conditions hold:
\begin{enumerate}
\item $\overline r$ is represented by a Fontaine--Laffaille module of weight $\lambda+\eta$, and the associated Serre weight $F(\lambda)$ is $5n$-generic: for every $j$,
\[
 \lambda_{j,1}-\lambda_{j,n}\leq p-6n,
 \qquad \lambda_{j,i}-\lambda_{j,i+1}\geq 5n
 \quad(1\leq i<n).
\]
In particular, the Hodge type lies in the strict Fontaine--Laffaille range.
\item $F(\lambda)$ is an extremal modular Serre weight for $\overline r$; equivalently, with the determinant twist used in \cite{LLHMPQ},
\[
 M_\infty\bigl(F(\lambda)\otimes_{\mathbb F}
 \omega^{n-1}\circ\det\bigr)\neq 0.
\]
\item For each $\tau_i$, the Kisin shape of $\overline r$ relative to $\tau_i$ is maximal:
\[
 \widetilde w^{*}(\overline r,\tau_i)=t_\eta
\]
in the notation of \cite{LLHMPQ}.
\end{enumerate}
\end{definition}

The last condition is a condition on the residual Galois representation, not on a characteristic-zero point. It can be checked from the matrices of the associated Fontaine--Laffaille module. The modularity condition is already part of the global patching datum. The numerical condition is stronger than merely being in the strict Fontaine--Laffaille range; this extra genericity is what makes the required weight elimination theorem available in arbitrary dimension.

\begin{thm}\label{thm:FL-evaluation}
Assume that $\overline r$ satisfies Definition \ref{def:GFL}. Then $\phi^\circ$ is surjective. Consequently $\phi$ is surjective and
\[
 X=\Sp R_\infty(\sigma_{max})[1/p].
\]
In particular, $X$ is Zariski closed.
\end{thm}

\begin{proof}
The patched functor on finite $\mathcal O[K]$-modules is exact, since $M_\infty$ is projective in the category of pseudocompact $\mathcal O[[K]]$-modules. By \cite[Lemma 10.2.8]{LLHMPQ}, the modularity of the $5n$-generic weight implies the required modularity of all obvious weights. For every rank-one factor above, the maximal-shape assumption and \cite[Lemma 10.2.12]{LLHMPQ} imply that the Jordan--H\"older factors of the two comparison cokernels are disjoint from the modular Serre weights of $\overline r$. Exactness therefore kills both cokernels after applying the patched functor.

It follows from \cite[Proposition 10.2.11]{LLHMPQ}, applied as in \cite[Proposition 10.2.13]{LLHMPQ}, that every normalized rank-one operator acts by a unit on the corresponding patched module. Here regularity of the auxiliary maximal-shape deformation rings is part of the cited argument; the crystalline deformation ring in the Fontaine--Laffaille range is formally smooth by \cite[Theorem A]{BooherLevin}. Thus the reduction of each factor, and hence the reduction of $\phi^\circ$, is surjective. The cokernel of $\phi^\circ$ is a finite module over the complete local patched deformation ring, so Nakayama's lemma gives $\operatorname{coker}(\phi^\circ)=0$. Inverting $p$ proves that $\phi$ is surjective, and the assertion about $X$ follows from its definition.
\end{proof}

\subsection{Existence of a $G$-invariant norm}\label{A.2.2}
Recall that we have fixed a residual Galois representation $\overline{r} : G_F \longrightarrow GL_n(\mathbb{F})$. The patching construction carried out in \cite{MR3529394} associates to $\overline{r}$ a modules $M_{\infty}(\overline{r}):=M_{\infty}$ and also the ring $R_{\infty}(\sigma_{min})(\overline{r}):=R_{\infty}(\sigma_{min})$.

Let $\rho : G_F \longrightarrow GL_n(E)$ be a potentially semi-stable Galois representation of weight $\sigma_{alg}$ and inertial type $\tau$. By the theory of Fontaine, the Galois representation $\rho$ corresponds to a filtered admissible $(\varphi,N, \mathrm{Gal}(L/F))$-module $\tilde{D}$. Proposition \ref{2.5} says that there is an admissible filtered $\varphi$-module $D$ such that the underlying Weil representation and the Hodge-Tate weights of $D$ are the same as that of $\tilde{D}$, but the monodromy operator is equal to zero. Let $r$ be the Galois representation corresponding to $D$. This Galois representation is potentially crystalline. Observe that the Galois representations $r$ and $\rho$ may have different reduction mod $p$ since the filtrations  on the corresponding $(\varphi,N)$-modules may differ. 

Let now $\overline{\rho} : G_F \longrightarrow GL_n(\mathbb{F})$ be the reduction mod $p$ of $\rho$. We can consider a new deformation problem associated to $\overline{\rho}$. Let $\tilde{M}_{\infty}:=M_{\infty}(\overline{\rho})$,   $\tilde{R}_{\infty}(\sigma_{min}):=R_{\infty}(\sigma_{min})(\overline{\rho})$ and let $\tilde{R}_{\infty}(\sigma_{min})'=R_{\infty}\otimes_{R_{\tilde{\mathfrak{p}}}^{\square}}R^{\square}_{\overline{\rho}}(\tau, \mathbf{v})$, where the ring $R^{\square}_{\overline{\rho}}(\tau, \mathbf{v})$ parametrizes all the potentially semi-stable lifts of $\overline{\rho}$ of weight $\sigma_{alg}$ and inertial type $\tau$.

Recall that we say that $r$ is generic when $\pi_{sm}(r)$ is generic.








\begin{thm}\label{3.4}
Let $r$ and $\rho$ be two Galois representations, as above. Assume that $F/\Qp$ is unramified, that the inertial type is trivial, and that the Hodge type is associated with a $5n$-generic Fontaine--Laffaille weight. Assume moreover that $\overline r$ satisfies Definition \ref{def:GFL}. Assume that $\rho$ is generic (we do not assume that $\rho$ corresponds to a point lying on an automorphic component) and that $r_x:=r$ corresponds to a closed point $x \in \Sp R_{\infty}(\sigma_{max})[1/p]$ (deformation problem for $\overline{r}$). Then $BS(r_x)$ and $BS(\rho)$ both admit a $G$-invariant norm. The completions of $BS(r)$ and $BS(\rho)$ with respect to these norms are admissible.
\end{thm}

\begin{proof} Recall that we have $V(r_x)^{l.alg}=\pi_x\otimes \pi_{alg}(r_x)$ by Proposition 4.33 \cite{MR3529394}, where $\pi_x$ is some smooth admissible representation in the Bernstein component $\Omega$. By Corollary 3.11 \cite{MR3529394}, the irreducible representation $\pi_{sm}(\rho)$ is the socle of $\cI_K^G \sigma_{max}(\lambda) \otimes_{\mathfrak{Z}_{\Omega},\chi_{\pi_{sm}(\rho)}} E$.

Let $y$ be the image of $x$ by the map $\Sp R_{\infty}(\sigma_{max})[1/p] \longrightarrow \Sp \mathfrak{Z}_{\Omega}$ (map $\alpha^{\sharp}$ of \cite[Theorem 5.6]{Pyv3} or equivalently the map from Theorem 4.19 \cite{MR3529394}). According to \cite[Theorem 3.3]{Pyv3} and \cite[Theorem 5.6]{Pyv3}, the closed point $y$ depends only on the eigenvalues of the linearised Frobenius $\varphi^f$ (which acts on both $D$ and $\tilde{D}$). The Galois representation $\rho$ corresponds to a closed point $\tilde{x}$ of $\Spm R^{\square}_{\overline{\rho}}(\tau, \mathbf{v})[1/p]$. Let $\tilde{y}$ be the image of $\tilde{x}$ by the map $\Sp R^{\square}_{\overline{\rho}}(\tau, \mathbf{v})[1/p] \longrightarrow \Sp \mathfrak{Z}_{\Omega}$(\cite[Theorem 3.3]{Pyv3}). The same conclusions hold for $\tilde{y}$.

By construction, $D$ and $\tilde{D}$ have the same action of $\varphi$. Then together with the observation above, it follows that $y=\tilde{y}$.

Let $\gamma_x:=\cI_K^G \sigma_{max}(\lambda) \otimes_{\mathfrak{Z}_{\Omega}} \kappa(y)$, this is a parabolic induction of a supercuspidal representation of a Levi subgroup of $G$. Let $\delta_x:= \cI_K^G \sigma_{min}(\lambda) \otimes_{\mathfrak{Z}_{\Omega}} \kappa(y)$. By definition $BS(r_x)= \gamma_x \otimes \pi_{alg}(r_x)$. 

Since $y=\tilde{y}$, it follows that $\cI_K^G \sigma_{max}(\lambda) \otimes_{\mathfrak{Z}_{\Omega},\chi_{\pi_{sm}(\rho)}} E = \gamma_x$. We conclude that $\pi_{sm}(\rho) = \mathrm{soc}(\gamma_x)$ is the socle of $\gamma_x$, and this irreducible representation is generic. Hence by \cite[Theorem 2.1]{Pyv2}, $\mathrm{Hom}_K(\sigma_{min}(\lambda), \mathrm{soc}(\gamma_x))\neq 0$. It follows that there is a non zero $G$-equivariant map $\delta_x \longrightarrow \mathrm{soc}(\gamma_x)$. Hence we have a non-zero map $\delta_x \longrightarrow \mathrm{soc}(\gamma_x) \hookrightarrow \gamma_x$.

Since $x \in \mathrm{Supp} (M_\infty(\sigma_{max}^{\circ}))$, because $r_x$ is potentially crystalline, we have $\mathrm{Hom}_K(\sigma_{max}(\lambda), \pi_x) \neq 0$, therefore there is a non-zero map $\gamma_x \longrightarrow \pi_x$. By Lemma \ref{3.2} and Lemma \ref{3.3} the composition $\delta_x \longrightarrow \mathrm{soc}(\gamma_x) \longrightarrow \gamma_x \longrightarrow \pi_x$ is non-zero if and only if the map $\gamma_x \longrightarrow \pi_x$ is injective. We will prove that $\gamma_x \longrightarrow \pi_x$ is injective.

With the notation introduced in the previous section, Theorem \ref{thm:FL-evaluation} gives $X=\Sp R_\infty(\sigma_{max})[1/p]$. Thus $x\in X$, and Lemma \ref{3.3} shows that $\gamma_x \longrightarrow \pi_x$ is injective.

So we have that $BS(r_x) \hookrightarrow V(r_x)^{l.alg}$, by a similar argument. The restriction of the norm on Banach space representation $V(r_x)$ induces a $G$-invariant norm on $BS(r_x)$. Since $\pi_{sm}(\rho)$ is the socle of $\gamma_x$, we have a $G$-equivariant injection $BS(\rho) \hookrightarrow BS(r_x)$. So we also obtain a $G$-invariant norm on $BS(\rho)$ by restricting a $G$-invariant norm on $BS(r_x)$.
\end{proof}

\noindent \textbf{Remark.} It is expected that $BS(r_x) \simeq V(r_x)^{l.alg}$. 

\subsection{Example}\label{A.2.3}

In this section, we give an example illustrating Theorem \ref{3.4}. Let $\rho$ and $r$ be as in that theorem. Under the assumptions of \cite[Corollary 5.5(2)]{MR3529394}, all components are automorphic. This removes the automorphic-component condition, while the residual Fontaine--Laffaille hypothesis of Definition \ref{def:GFL} is still required. We now specify $r$ and $\rho$.

Let $F=\Qp$, $n=3$ and  let $r$, $s$ be two integers such that $0<r<s$, $r\leq p-1$ and $s-r \leq p-1$. Let $v_p$ be a valuation $\overline{\Q}_p$ with $v_p(p)=1$. Assume $p\neq 3$. Let $\tilde{D}$ be an admissible filtered $(\varphi, N)$-module with Hodge-Tate weights $0<r<s$, from Example 3.40 \cite{MR3459578}: 

\begin{enumerate}
\item[$\bullet$] $\mathrm{Fil}^{r} \tilde{D} = E(e_1,e_2)$ and $\mathrm{Fil}^{s} \tilde{D} = E(e_1)$.
\item[$\bullet$] $N=\begin{pmatrix}
0 & 0 & 0\\ 
0 & 0 & 0\\
1 & 0 & 0
\end{pmatrix}$ and $\varphi =\begin{pmatrix}
p\mu & 0 & 0\\ 
0 & \mu_2 & 0\\
0 & 0 & \mu
\end{pmatrix}$, with $\mu \neq \mu_2\neq p^{\pm 1}\mu$.
\item[$\bullet$] $ r-1 > v_p(\mu) > (r-1)/2 >0$ and $1+v_p(\mu) \geq v_p(\mu_2) \geq v_p(\mu)$.
\end{enumerate}

According to Proposition 3.41 \cite{MR3459578}, the $(\varphi, N)$-module above is irreducible, because $r-1 > v_p(\mu) > (r-1)/2$ and $r>1$. 

Let $\rho$ be a semi-stable Galois representation with Hodge-Tate weights $0<r<s$ corresponding to $\tilde{D}$. Let $P$ be a standard parabolic subgroup of $G$ corresponding to a partition $(2,1)$, let $\chi$ and $\chi_2$ unramified characters of $\Qp$ such that $\chi(p)=\mu$ and $\chi_2(p)=\mu_2$. Then 
$$\pi_{sm}(\rho) = i_P^G((St\otimes |.|^{-1/2}\chi)\otimes \chi_2)\otimes|\det|,$$

\noindent where $St$ is the Steinberg representation of $GL_2(\Qp)$.  Since after killing the monodromy there is not a unique choice of a filtration that makes the underlying $\varphi$-module admissible, we may choose $r$, so that the Galois representation $r$ corresponds to a $\varphi$-module $D$, from Example 2.61 \cite{MR3459578}:

\begin{enumerate}
\item[$\bullet$] $\mathrm{Fil}^{r} D = E(e_1+e_2+e_3, e_2+2e_3)$ and $\mathrm{Fil}^{s} D = E(e_1+e_2+e_3)$.
\item[$\bullet$] $N=0$ and $\varphi =\begin{pmatrix}
p\mu & 0 & 0\\ 
0 & \mu_2 & 0\\
0 & 0 & \mu
\end{pmatrix}$, with $\mu \neq \mu_2\neq p^{\pm 1}\mu$.
\item[$\bullet$] $s>r>1+v_p(\mu) \geq v_p(\mu_2) \geq v_p(\mu)>(r-1)/2>0$ and $1+v_p(\mu) + v_p(\mu_2) + v_p(\mu) = r+s$. 
\end{enumerate}

Since $s>1+v_p(\mu) \geq v_p(\mu_2) \geq v_p(\mu)>0$, it follows from Proposition 2.62 \cite{MR3459578} that this $\varphi$-module is irreducible. Then by classical Langlands correspondence: 
$$\pi_{sm}(r) = L(|.|^{-1}\chi,\chi_2,\chi)\otimes |\det| = L(\chi, |.|\chi_2, |.|\chi)$$

\noindent  is a Langlands quotient where the segments $\chi$ and $|.|\chi$ are linked. 

\begin{prop}\label{3.10}
Let $\rho$ be a generic semi-stable Galois representation and $r$ be a crystalline Galois representation as above, both having same Hodge-Tate weights $0<r<s$ and the same action of the Frobenius $\varphi$. Assume that the Serre weight associated with these Hodge--Tate weights is $15$-generic Fontaine--Laffaille and that the reduction $\overline r$ of a $G_{\Qp}$-stable lattice in $r$ satisfies Definition \ref{def:GFL}. Let $\chi$ and $\chi_2$ be two unramified characters such that $\chi(p) = \mu$ and $\chi_2(p) = \mu_2$, with $\mu$ and $\mu_2$ satisfying relations:
\begin{enumerate}
\item[(A)] $\mu \neq \mu_2\neq p^{\pm 1}\mu$ and $r>1+v_p(\mu) \geq v_p(\mu_2) \geq v_p(\mu)>(r-1)/2>0$.
\item[(B)] $1+v_p(\mu) + v_p(\mu_2) + v_p(\mu) = r+s$.
\end{enumerate}
Then $BS(r):=(i_B^G(\chi\otimes|.|\chi_2 \otimes |.|\chi)) \otimes \pi_{alg}(r)$ admits a $G$-invariant norm, where $B$ is a Borel subgroup of $G$ and $\pi_{alg}(r)$ is an irreducible algebraic representation of highest weight $\psi(\mathrm{diag}(t_1,t_2,t_3)) = t_1^{0}.t_2^{-r+1}.t_3^{-s+2}$ with respect to upper triangular Borel $B$. Moreover, $BS(\rho):=\pi_{sm}(\rho) \otimes \pi_{alg}(r)$ admits also a $G$-invariant norm and the completions of $BS(r)$ and $BS(\rho)$ with respect to these norms are admissible.
\end{prop}

\begin{proof} First observe that $\pi_{sm}(r)$ is an irreducible quotient of $i_B^G(\chi\otimes|.|\chi_2 \otimes |.|\chi)$, hence by definition $BS(r):=(i_B^G(\chi\otimes|.|\chi_2 \otimes |.|\chi)) \otimes \pi_{alg}(r)$. By \cite[Corollary 5.5(2)]{MR3529394}, the point defined by $r$ lies on an automorphic component. The additional assumptions in the statement are precisely the residual hypotheses needed in Theorem \ref{3.4}; hence that theorem gives a $G$-invariant norm on $BS(r)$. Since $BS(\rho) \hookrightarrow BS(r)$, restriction gives a $G$-invariant norm on $BS(\rho)$ as well.
\end{proof}

We will prove that the $G$-invariant norm on $BS(r)$ does not come from a restriction of a norm on a parabolic induction of a unitary character. First, we will prove the following lemma:

\begin{lemma}\label{3.6}
Let $\pi$ be a smooth admissible representation of $G$ and $\sigma$ an algebraic irreducible representation of $G$ with highest weight $\psi$ with respect to upper triangular Borel subgroup $B$ of $G$. Then $(\pi \otimes \sigma)_N \simeq \pi_N \otimes\sigma_N$.
\end{lemma}

\begin{proof} Let $V$ denote vector space equipped with a $G$-action. We will denote by $V(N)$ the space spanned by $n.v-v$, $n \in N$ and $v \in V$, and by $V_N=V/V(N)$. We will identify injective maps with the inclusions. Since $(\pi \otimes \sigma)(N) \subseteq \pi(N)\otimes \sigma(N)$ then we have $(\pi \otimes \sigma)_N \twoheadrightarrow \pi_N \otimes\sigma_N$.

The representation $\sigma$ is finite dimensional. Let $w$ be the highest weight vector of $\sigma$. Observe that $\sigma_N$ is one dimensional generated by $w$. 

Let $v \in \pi$. Since $\pi$ is smooth, the vector $v$ is fixed by some compact open $N_0 \subseteq N$. We have also that $\sigma_N=\sigma_{N_0}=E.w$, because this representation is algebraic. Since $\sigma$ is finite dimensional, we may choose $w_1,\ldots w_d$ a basis of $\sigma$ such that $w_m \in \sigma(N_0)$ for $m\neq d$ and $w_d=w$. Then $w_m \in \sigma(N_0)$ can be written as $w_m=\sum\limits_{k=1}^{d} a_k (n_k-1) w_k$, where $a_k$ are some scalars and $n_k\in N_0$. It follows that for $m\neq d$:
$$v\otimes w_m = v\otimes \sum\limits_{k=1}^{d} a_k (n_k-1) w_k = \sum\limits_{k=1}^{d} a_k (v\otimes(n_k-1) w_k)=$$ 
$$= \sum\limits_{k=1}^{d} a_k ((v\otimes n_k w_k)- v\otimes w_k)= \sum\limits_{k=1}^{d} a_k ((n_k v\otimes n_k w_k)- v\otimes w_k)=$$
$$ =\sum\limits_{k=1}^{d} a_k (n_k(v\otimes w_k)-v\otimes w_k)= \sum\limits_{k=1}^{d} a_k (n_k-1)(v\otimes w_k) \in (\pi \otimes \sigma) (N)$$

\noindent This shows that $\pi\otimes(\sigma(N)) \subseteq (\pi \otimes \sigma)(N)$. Therefore we get a surjection $\pi \otimes \sigma_N \twoheadrightarrow (\pi \otimes\sigma)_N$. Since $N$ acts trivially on $\sigma_N$, this map factors through $\pi_N \otimes \sigma_N$, as 
$$\xymatrix{
\pi \otimes \sigma_N \ar@{->>}[d] \ar@{->>}[r] & (\pi \otimes \sigma)_N \\
\pi_N \otimes \sigma_N \ar@{->>}[ur] & }$$
\noindent The composition map $\pi_N \otimes \sigma_N\twoheadrightarrow(\pi \otimes \sigma)_N \twoheadrightarrow \pi_N \otimes \sigma_N$ is the identity. This allows us to conclude.\end{proof}

\begin{prop}\label{3.11}
Let $T$ be a group of diagonal matrices and $N$ group of unipotent matrices such that both are subgroups of a Borel $B=TN$ as in Proposition \ref{3.6}. Let $\theta : T \rightarrow \mathcal{O}^{\times} \rightarrow E^{\times}$ be a unitary character. Then there is no such unitary character $\theta$, such that the we have an embedding $BS(r) \hookrightarrow \mathrm{Ind}_{B}^{G}(\theta)_{cont}$ (index $cont$ means that we consider continuous functions in this space). We also have that for any unitary character $\theta$ as above, there is no injection  $BS(\rho) \hookrightarrow  \mathrm{Ind}_{B}^{G}(\theta)_{cont}$.
\end{prop}

\begin{proof} Since $BS(r)$ is locally algebraic, we can restrict ourself to the space $ \mathrm{Ind}_{B}^{G}(\theta)^{l.an}$(the functions in this space are locally analytic, see \cite{MR2392361} for definitions), and work in the category of locally analytic representations. According to Theorem 4.2.6 \cite{MR1691735}, we have a Frobenius reciprocity in the category of locally analytic representations:
$$\mathrm{Hom}_{G}^{cont}(BS(r), \mathrm{Ind}_{B}^{G}(\theta)^{l.an}) \simeq \mathrm{Hom}_{B}^{cont}(BS(r)|B, \theta)$$

Let $\delta^{1/2}=|.|^{-1/2}\otimes 1 \otimes |.|^{1/2}$ be a square root of the modulus character. Since $\theta$ is trivial on $N$, $BS(r)$ factors through $N$-coinvariants, then:
$$\mathrm{Hom}_{G}^{cont}(BS(r), \mathrm{Ind}_{B}^{G}(\theta)^{l.an}) \simeq \mathrm{Hom}_{T}^{cont}((BS(r))_N, \theta)$$

\noindent Then Lemma \ref{3.6} allows us to compute:
$$(BS(r))_N = \delta^{1/2}.r_B^G(i_B^G(\chi\otimes |.|\chi_2 \otimes |.|\chi)) \otimes \pi_{alg}(r)_N,$$

\noindent where $r_B^G$ is the left adjoint functor of $i_B^G$ in the category of smooth representations. Let $\tilde{\chi}:= \chi\otimes|.|\chi_2 \otimes |.|\chi$. Theorem 1.2 \cite{MR584084} tells us that there is a filtration:
$$0=\tau_0 \subset \tau_1 \subset ... \subset \tau_6 = r_B^G(i_B^G(\tilde{\chi})),$$

\noindent such that $\tau_i/\tau_{i-1} \simeq \tilde{\chi}^{w_i}$, where $w_i$ is an element of the symmetric group in 3 letters $\mathfrak{S}_3= \left\lbrace w_1,...,w_6 \right\rbrace$. Let $\psi(\mathrm{diag}(t_1,t_2,t_3)) = t_1^{0}.t_2^{-r+1}.t_3^{-s+2}$ be the weight of the highest weight representation $\pi_{alg}(r)$ with respect to upper triangular matrices. If we have 
$$\mathrm{Hom}_{T}^{cont}(\delta^{1/2}.\tilde{\chi}^w.\psi, \theta) = 0,$$

\noindent for all $w \in \mathfrak{S}_3$, then $\mathrm{Hom}_{G}^{cont}(BS(r), \mathrm{Ind}_{B}^{G}(\theta)^{l.an})=0$. It is enough to prove that for any $w \in \mathfrak{S}_3$ the character $\delta^{1/2}.\tilde{\chi}^w.\psi$ is not unitary. Indeed we have that:
$$v:=v_p(\delta^{1/2}.\tilde{\chi}^w.\psi(\mathrm{diag}(1,p,1)))\in$$ 
$$\left\lbrace 1-r+v_p(\mu), -r + v_p(\mu),-r+v_p(\mu_2)\right\rbrace ,$$

\noindent and in all these cases we have $v<0$ by relations (A) of Proposition \ref{3.10}. It follows that the character $\delta^{1/2}.\tilde{\chi}^w.\psi$ does not take its values in $\mathcal{O}^{\times}$. For $BS(\rho)$ we repeat the same proof.
\end{proof}

The locally algebraic representations $BS(\rho)$ and $BS(r)$ have the same central character. We can easily verify that, by relation (B) of Proposition \ref{3.10}, we have $v_p(\tilde{\chi}.\psi(\mathrm{diag}(p,p,p)))= v_p(\mu)+v_p(\mu_2)-1+v_p(\mu)+3-r-s=2v_p(\mu)+v_p(\mu_2)+1-r-s=0$, therefore the central character of $BS(r)$ (or of $BS(\rho)$) is unitary.

We have shown that, in this example, $BS(r)$ can not be embedded into $\mathrm{Ind}_{B}^{G}(\theta)_{cont}$, with $\theta$ unitary. Thus the $G$-invariant norm on $BS(r)$, obtained by the Theorem \ref{3.4}, does not come from a restriction of a $G$-invariant norm on a parabolic induction of a unitary character. Same conclusions hold for $BS(\rho)$.

\subsection*{Acknowledgments}  The results of this paper are an improvement of the last part of the author's PhD thesis. The author tremendously grateful to his advisor Vytautas Pa\v{s}k\={u}nas for sharing his ideas with him, and also the author would like to thank Qijun Yan for some productive discussions. This work was funded by Morningside Center of Mathematics, CAS, and partially supported by SFB/TR 45 of the DFG.

\bibliographystyle{alpha}
\addcontentsline{toc}{section}{References}
\bibliography{Locally_alg_cryst_non_gen}
\nocite{*}

\textit{E-mail address}: a413xpyv@hotmail.com

\end{document}